\documentclass[12pt]{article}
\usepackage{amssymb}
\usepackage{amsmath,amsfonts,amsthm}
\usepackage{color}

\newcommand{\ba}{$$\begin{array}}\newcommand{\ea}{\end{array}$$}
\newcommand{\bea}{\begin{eqnarray*}}\newcommand{\eea}{\end{eqnarray*}}
\newcommand{\be}{\begin{equation}}\newcommand{\ee}{\end{equation}}
\newcommand{\bd}{\begin{definition}} \newcommand{\ed}{\end{definition}}
\newcommand{\brs}{\begin{remarks}\rm}   \newcommand{\ers}{\end{remarks}}
\newcommand{\br}{\begin{remark}\rm}     \newcommand{\er}{\end{remark}}
\newcommand{\bt}{\begin{theorem}}       \newcommand{\et}{\end{theorem}}
\newcommand{\bl}{\begin{lemma}}         \newcommand{\el}{\end{lemma}}
\newcommand{\bco}{\begin{corollary}}    \newcommand{\eco}{\end{corollary}}
\newcommand{\bp}{\begin{proposition}}   \newcommand{\ep}{\end{proposition}}
\newcommand{\bo}{\begin{observation}\rm}\newcommand{\eo}{\end{observation}}
\newcommand{\bex}{\begin{examples}\rm}   \newcommand{\eex}{\end{examples}}
\newcommand{\bos}{\begin{observations}\rm}\newcommand{\eos}{\end{observations}}
\newcommand{\bx}{\begin{example}\rm}   \newcommand{\ex}{\end{example}}
\newcommand{\bexe}{\begin{exercise}\rm}   \newcommand{\eexe}{\end{exercise}}
\newcommand{\bpf}{\begin{proof}\rm}   \newcommand{\epf}{\end{proof}}

\newtheorem{definition}{Definition}[section]
\newtheorem{theorem}[definition]{Theorem}
\newtheorem{lemma}[definition]{Lemma}
\newtheorem{corollary}[definition]{Corollary}
\newtheorem{proposition}[definition]{Proposition}
\newtheorem{remarks}[definition]{Remarks}
\newtheorem{remark}[definition]{Remark}
\newtheorem{observation}[definition]{Observation}
\newtheorem{examples}[definition]{Examples}
\newtheorem{observations}[definition]{Observations}
\newtheorem{example}[definition]{Example}
\newtheorem{exercise}[definition]{Exercise}

\begin{document}
\begin{center} {\bf Some results on $K$-algebras} \end{center}
\begin{center}Pramod K. Sharma\\ e--mail:
 pksharma1944@yahoo.com\\
 School of Mathematics, Vigyan Bhawan, Khandwa Road,\\ INDORE--452
017, INDIA.\end{center} \section*{Abstract:}We give a new proof of
the classical result due to Rodney Y. Sharp and Peter Vamos on the
dimension of tensor product of a finite number of field extensions
of a given field.\section{Introduction} Let $K$ be a field. In this
note, we prove some results on $K$-algebras. All rings and algebras
are commutative with identity $\neq 0$. By the dimension of a ring
$A$ we mean the Krull dimension and denote it by dim $A$. The
transcendence degree of a field extension $L/K$ shall be denoted by
$ {tr deg_K L}$. The results in this note grew while trying to
understand the classical result on dimension of the tensor product
of two field extensions proved in [6]. We first prove [Theorem 1] :
 Let $R\subset A$ be rings where $R$ is an integral domain with its field of fraction $K$. Then
(1) If $X_1,X_2,...,X_n$ are algebraically independent over $A$ and
 $A$ contains $t_1, t_2, ...t_n$ algebraically independent over $R$
then for $L = K (X_1,..., X_n), \,\, dim (L\bigotimes_R A) \geq n +
dim \,\,   S^{-1} A$  where $S$ is the multiplicatively closed
subset $R[t_1,...,t_n]-\{0\}$ of $A$, and (2) If
$X_1,X_2,...,X_n,...$ are algebraically independent over $A$ and $A$
contains $t_1, t_2, ...t_n,...$ algebraically independent over $R$
then for $L=K(X_1,..., X_n,...),\,\,  dim (L\bigotimes_R A)=
\infty$. In Corollary 2.3, it is shown that equality holds in
Theorem 1 under certain conditions. These results are used to find
the dimension of the tensor product of a finite number of field
extensions of a given field proved in [7]. Further, we give [Theorem
2.7] an alternative proof of the well known result that for an
affine $K$-algebra $A$ over a field $K$, for any non-zero-divisor
$f\in A$, $dim A = dim A[1/f]$.

\section { Main Results}

Before we prove that main results, let us recollect :\\ (i) [5,
Theorems 7.3 and 9.5]: If $B$ is a faithfully flat $A$-algebra
then $dim B\geq dim A.$\\ (ii) [5, Exercise 9.2] If a ring $B$ is
an integral extension of a ring $A$ then $dim A=dim B$. \\

We shall use these facts, whenever required, without further
mention.
 \bt Let $R\subset A$ be rings where $R$ is an
integral domain. Let $K$ be the field of fractions of R. Then \\(1)
If $X_1, \cdots , X_n$
 are algebraically independent over $A$ and $A$ contains $t_i,$ $i=1, \cdots , n$
  algebraically independent over $R$, then
  $$\dim K(X_1, \cdots , X_n)\otimes_RA
\geq n+\dim S^{-1}A $$ where $S=R[t_1, \cdots , t_n]-\{0\}$.
Further, if $A$ is Noetherian, then $$\dim K(X_1, \cdots ,
X_n)\otimes_R A \leq dim A+n. $$ (2) If $X_1, \cdots , X_n,\cdots$
are algebraically independent over $A$ and $A$ contains $t_i,
i=1,2,\cdots,n\cdots$ algebraically independent over $R$, then
$$\dim K(X_1, \cdots , X_n, \cdots)\otimes_R A=\infty. $$ \et \bpf
(1) Let $P_0'\subsetneqq P_1'\subsetneqq P_2'\subsetneqq \cdots
\subsetneqq P_m'$ be a chain of prime ideals in $S^{-1}A$. Then
there exist prime ideals $P_0\subsetneqq P_1\subsetneqq
P_2\subsetneqq \cdots \subsetneqq P_m$ in $A$ such that $P_i\cap S=
\phi$ and $S^{-1}P_i =P_i'$. Note that $$P_0\subsetneqq
P_1\subsetneqq P_2\subsetneqq \cdots \subsetneqq P_m\subsetneqq
(P_m,\, X_1-t_1)\subsetneqq \cdots \subsetneqq (P_m, X_1-t_1,
\cdots, X_n-t_n)$$ is a chain of prime ideals in $A[X_1, \cdots ,
X_n]$. If for $T=R[X_1, \cdots , X_n]-\{0\}$, $T\cap(P_m, X_1-t_1,
\cdots, X_n-t_n)\neq \phi$, then there exist $f(X_1, \cdots ,
X_n)(\neq 0)\in R[X_1, \cdots , X_n]$ such that $$f(X_1, \cdots ,
X_n)=g(X_1, \cdots , X_n)+\sum (X_i-t_i)h_i(X_1,\cdots , X_n)$$
where $h_i\in A[X_1, \cdots , X_n]$ and $g(X_1, \cdots , X_n)\in
P_m[X_1, \cdots , X_n]$. This implies that $f(t_1, \cdots ,
t_n)=g(t_1, \cdots , t_n)\in P_m $. Since $t_i$'s are algebraically
independent over $R$, $f(t_1, \cdots , t_n)\neq 0\in P_m\cap S$.
This contradicts  our assumption on $P_i$'s. Therefore $ T\cap (P_m,
X_1-t_1, \cdots , X_n-t_n)=\phi$, and
$$\dim T^{-1}(A[X_1, \cdots , X_n])\geq n+\dim S^{-1}A $$ where
$T=R[X_1, \cdots , X_n]-\{0\}$. Now, note that $$R[X_1, \cdots ,
X_n]\otimes_RA\cong A[X_1, \cdots , X_n]$$ as $R[X_1, \cdots ,
X_n]$-algebras.  Hence  $$K(X_1, \cdots , X_n) \otimes_{R[X_1,
\cdots , X_n]}A[X_1, \cdots , X_n] \cong T^{-1}A[X_1, \cdots,
X_n]$$ $$\Rightarrow \,\,\dim(K(X_1, \cdots , X_n)\otimes_RA)\geq
n+\dim S^{-1}A.$$

The final part of the statement is immediate since $K(X_1,...,X_n
) \otimes _R A$ is a localization of $R[X_1, \cdots ,
X_n]\otimes_RA$ which is isomorphic to $A[X_1, \cdots , X_n].$
 \\
Further, as $A$ is Noetherian, $ dim A[X_1, \cdots , X_n]= dim A+n
  $ [5, Theorem 15.4]  \\  (2) Let us note that $$K(X_1, \cdots ,
X_n,\cdots) \otimes_{K(X_1, \cdots , X_n)}(K(X_1, \cdots , X_n)
\otimes_R A) \cong K(X_1, \cdots , X_n, \cdots)\otimes_R A$$ Hence
$ K(X_1, \cdots , X_n,\cdots )\otimes_R A$ is faithfully flat
$K(X_1, \cdots , X_n  )\otimes_R A$ -  \,\,\, algebra.
\\ Therefore $$dim K(X_1, \cdots , X_n, \cdots)\otimes_R A\geq dim
K(X_1, \cdots , X_n )\otimes_R A$$ $$\,\,\,\,\,\,\,\,\,\,\,\,\,\,
\,\,\,\,\,\,\, \,\,\,\,   \, \geq n \,\,\,\,\,\,\,\,\,\,\,
(use(1))$$ $$ \Rightarrow dim K(X_1, \cdots , X_n,
\cdots)\otimes_R A = \infty.$$ \epf

 \br
In above Theorem, if $B$ is any $K(X_1, \cdots , X_n)$-algebra,
then $$  dim B \otimes_R A \geq dim K(X_1, \cdots , X_n
 )\otimes_R A$$ $$\geq n+\dim S^{-1}A $$ Further, if $B$ is
$K (X_1, \cdots , X_n, \cdots )-$ algebra, then $$dim B \otimes_R
A = \infty.$$

These observations are immediate since $B  \otimes_R    A$ is
faithfully flat $K(X_1, \cdots , X_n )\otimes_R A (K(X_1, \cdots ,
X_n, \cdots)\otimes_R A)- algebra$. \er

\bco Let $K$ be a field and  $A$ be a  $K $-algebra. If $X_1, \cdots
, X_n $ are algebraically independent over $A$ and $A$ contains a
field extension of $K$ of transcendental degree $\geq n$, then $$dim
K(X_1, \cdots , X_n )\otimes_K A \geq n + dim A.$$ Further, if A is
Noetherian, then $$dim K(X_1, \cdots , X_n
 )\otimes_K A =n + dim A.$$ \eco

 \bpf  By assumption on $A$, there exist $t_1, \cdots , t_n$
algebraically independent over $K$ such that $K (t_1, \cdots , t_n)
\subset A$. Hence for $S= K [t_1, \cdots , t_n] -{0}, S^{-1} A=A$.
Therefore, by the Theorem 1,$$dim K(X_1, \cdots , X_n
 )\otimes_K A \geq n + dim A.$$ Further, let $A$ be
Noetherian. Then as $$ K(X_1, \cdots , X_n)\otimes_K A \cong
T^{-1} A [X_1, \cdots , X_n]$$ where $T = [X_1, \cdots ,
X_n]-{0}$, it is immediate that $$ dim K(X_1, \cdots ,
X_n)\otimes_K A \leq dim A[X_1, \cdots , X_n]$$ $$= n + dim A$$
Consequently $$ n + dim A = dim K(X_1, \cdots , X_n) \otimes_K
A.$$  \epf

\bt Let $L_i, i=1,\cdots,n$ be a field extension of a given field
$K$ and let $trgdeg_KL_i = t_i$. Assume $t_1\leq t_2 \cdots\leq
t_n{_{-1}}\leq t_n$. If $t_n{_{-1}} <\infty$ then $$ dim
(L_1\otimes_K \cdots \otimes_K L_n) =t_i + t_2 + \cdots +
t_{n-1},$$ otherwise $$ dim(L_1 \otimes _K \cdots \otimes_K L_n) =
\infty.$$ \et  \bpf We shall consider the two cases separately. \\
Case 1.  $t_1 \leq t_2 \leq \cdots\leq t_{n-1}<\infty.$
\\ Let $B_k = \{x_{k_{1}},x_{k_{2}},\cdots x_{kt_k} \}$ be a transcendental
basis of $L_k$ over $K$ for $k=1,2,\cdots,n-1$.
 Put
$E_k=K(x_{k1},x_{k2},\cdots x_{kt _k} )$. Then $E_k/K$ is purely
transcendental field extension of  transcendental degree $t_k$ and
$L_k/E_k$ is algebraic. Hence  $$ E_1 \otimes_K E_2 \otimes_K
\cdots \otimes_K E_{n-1} \otimes_K L_n \,\,\,\, \stackrel {{i_1
\otimes \cdots \otimes i_{n-1} \otimes Id}}
{\hookrightarrow}\,\,\,\, L_1\otimes_k L_2\otimes _K \cdots
\otimes _K L_n,$$ where $i_k : E_k\hookrightarrow L_k$ is
inclusion map for $k=1,\cdots,n-1$ and $Id$ is identity map, is an
integral extension. Therefore $$ dim (L_1\otimes_K \cdots
\otimes_K L_n)=dim (E_1 \otimes_K E_2 \otimes_K \cdots \otimes_K
E_{n-1} \otimes_K L_n).$$ Let $Y_{11},Y_{12},\cdots Y_{1t_1},
Y_{21} \cdots Y_{2t_2}, \cdots, Y_{(n-1)1},\cdots,
Y_{(n-1)t_{(n-1)}}$ be algebraically independent elements over
$K$. Then for $F_k = K(Y_{11}, \cdots,Y_{1t_k}), k=1,\cdots,n-1$,
we have $$ F_1 \otimes_K \cdots \otimes_K F_{n-1} \otimes_{K} L_n
\cong E_1\otimes_K \cdots \otimes_K E_{n-1} \otimes_K L_n $$
Therefore $$ dim (F_1 \otimes_K \cdots \otimes_K F_{n-1}
\otimes_KL_n) = dim(L_1 \otimes_K \cdots \otimes_K L_n)$$ Let us
note that $F_2 \otimes_K \cdots\otimes_K F_{n-1} \otimes_KL_n$ is
a localization of
 $$L_n [Y_{21} \cdots Y_{2t_2},\cdots,Y_{n-1,1,}\cdots,
 Y_{(n-1)t_{(n-1)}}]$$
  over a multiplicatively closed subset, hence is a Noetherian ring.
Therefore by Corollary 2.3, $$ dim ( F_1 \otimes_K \cdots
\otimes_K F_{n-1} \otimes_KL_n)= t_1 + dim (F_2 \otimes_K \cdots
\otimes_K F_{n-1} \otimes_KL_n).$$ By successive application of
the Corollary 2.3 or by induction it is immediate that $$ dim F_2
\otimes_K \cdots \otimes F_{n-1} \otimes_KL_n=t_2+\cdots+t_{n-1}$$
Hence in this case the result follows.
\\Case 2.   $t_{n-1}=t_n=\infty.$
\\First of all, note that for any $\sigma \in S_n$
$$ L_1  \otimes_K \cdots \otimes_K  L_n \cong L_{\sigma(1)}
\otimes_K \cdots \otimes_K  L_{\sigma (n)}. $$ Therefore $$ L_1
\otimes_K \cdots \otimes_K  L_n \cong L_n \otimes_K L_{n-1}
\otimes_K \cdots \otimes_K L_2  \otimes_K L_1.$$ Put $B=L_{n-1}
\otimes_K \cdots \otimes_K L_2  \otimes_K L_1$. Then $$ dim (L_1
\otimes_K \cdots \otimes_K  L_n )=dim L_n \otimes _KB.$$ By
assumption $B$ contains infinite algebraically independent
elements over $K$. Hence the result is immediate from Theorem
1(2).\\

\epf

 \br If $A_i, i=1,\cdots,n$ denote integral extension of
$L_i$, then $$dim A_1\otimes_K\cdots\otimes_KA_n=dim
L_1\otimes_K\cdots \otimes_KL_n.$$ Further, if $A_i$ is any
$L_i-algebra,$ then $$dim A_1\otimes_K\cdots\otimes_KA_n\geq dim
(L_1\otimes_K\cdots \otimes_KL_n).$$  \bl  Let $K[X_1,\cdots,X_n]=
K\mathbf{[\underline{X}]}$ be a polynomial ring in $n$-variables
$X_i, i=1,\cdots n$ over a field $K$. Then for any $f(\neq 0)\in
K[\mathbf{\underline{X}}],dim K[\underline{\mathbf{X}},1/f]=n.$ \el
\er \bpf Let $\overline{K}$ be the algebraic closure of $K$. Then,
since $\overline{K}[\underline{\mathbf{X}},1/f]$is integral over
$K[\underline{\mathbf{X}},1/f]$, we have $$dim
\overline{K}[\underline{\mathbf{X}},1/f]=dim
K[\underline{\mathbf{\mathbf{X}}},1/f]. $$ Hence, to prove the
result, we can assume that $K$ is algebraically closed. Note that
$dim K[\underline{\mathbf{X}}]=n$ and for the multiplicatively
closed subset $ S= \{f^{t}|t\geq 0 \},
S^{-1}K[\underline{\mathbf{X}}]=K[\underline{\mathbf{X}},1/f]$.
Since $f\neq 0, f$ does not vanish on $K^n$. Thus, if for
$\underline{ \lambda}=\lambda_1,\cdots,\lambda_n $    in    $K^n,
f(\underline{\lambda})\neq 0,$ then for the maximal ideal
$M=(X_1-\lambda,\cdots,X_n-\lambda_n)$ in $K[\underline{X}], M \cap
S=\phi$. Therefore $S^{-1} M$ is a maximal ideal in
$S^{-1}K[\underline{X}]$. Clearly, height of $M, i.e.\,\,\,
     htM=n=htS^{-1}M$. Therefore $dim K[\underline{\mathbf{X}},1/f]=n$.
\epf

\bt Let $A$ be an affine algebra over a field $K$. Then for any non-
zero-divisor $f$ in $A$, $dim A[1/f] = dim A$.\et

\bpf Let $A = \frac{K[X_1,\cdots,X_n]}{I}$. Since $f$ is a non-
zero-divisor in $A$, $f$ lies in no prime ideal associated to $I$ in
$K[X_1,\cdots,X_n]$. Let $p$ be an associated prime ideal of $I$ in
$K[X_1,\cdots,X_n]$ such that $$ dim A=dim
\frac{K[X_1,\cdots,X_n]}{p}.$$ Then $\overline{f}$, image of $f$ in
$\frac{K[X_1,\cdots,X_n]}{p}$, is non-zero. Note that $dim A[/f]\leq
dim A$. Further, as $\frac{K[X_1,\cdots,X_n]}{p}. [1/\overline{f}]]$
is a  quotient ring  of $A [1/f]$ in a  natural way, $$dim
A[1/f]\geq dim\frac{K[X_1,\cdots,X_n]}{p}[1/\overline{f}].$$

Thus to prove Theorem, it is sufficient to show that $$ dim
\frac{K[X_1,\cdots,X_n]}{p}= dim \frac{K[X_1,\cdots,X_n]}{p}
[1/\overline{f}].$$

Let us observe that $$\theta :\frac{K[X_1,\cdots,X_n][Y]}{(p,f
Y-1)}\rightarrow \frac{K[X_1,\cdots,X_n]}{p} [1/\overline{f}]$$
$$Y \mapsto 1/\overline{f}$$ is $\frac{K[X_1,\cdots,X_n]}{p}$
algebra isomorphism. Therefore $$dim \frac{K[X_1,\cdots,X_n]}{p}
[1/\overline{f}]=dim \frac{K[X_1,\cdots,X_n][Y]}{(p,f Y-1)}.$$ We
note that $fY-1 \not\in p[Y]$. As A $\frac{K[X_1,\cdots,X_n]}{p}
[1/\overline{f}]$ is an integral domain, the ideal $(p,f Y-1)$ is
prime in $K[X_1,\cdots,X_n, Y]$. Now, note that $K[X_1,\cdots,X_n,
Y]$ is a Cohen-Macaulay ring of dimension $n+1$. By [4,Ex. 19,page
104], $ht (p,f Y-1)=ht p +1.$ Therefore $$ dim
\frac{K[X_1,\cdots,X_n][Y]}{(p,f Y-1)}=(n+1)-(ht p+1)$$ $$=n-ht
p$$ $$= dim \frac{K[X_1,\cdots,X_n]}{p}.$$ Thus $dim A=dim
A[1/f].$ \epf   We, now, deduce the following well known result:

 \bco  Let $A$ be an affine algebra over a field
$K$ which is an integral domain. Then $dim  A=trdeg_K L$ where $L$
is the field of fractions of $A$.\eco  \bpf  Let $\{y_1, \cdots,
y_s\}$ be a maximal algebraically independent set of elements in $A$
over $K$. Then every a $\in A$ is algebraic over $K[y_1, \cdots
y_s]$. Since $A$ is an affine algebra over $K, A=K[a_1, \cdots,
a_t]$ for some $a_i, i=1,2, \cdots, t$. Since each $a_i$ is
algebraic over $K[y_1, \cdots y_s]$ there exists an element $f(\neq
0) \in K[y_1, \cdots y_s]$ such that $A[1/f]$ is integral over
$K[y_1, \cdots y_s][1/f]$. Thus $$dim A[1/f]=dim K[y_1, \cdots,
y_s][1/f]$$
$$\,\,\,\,\,\,\,\,\,=s \,\,\,\,\,\,\,\,\,\,\,\,   (Lemma 2.6)$$
Therefore by  Theorem, it is immediate that $dim A=trdeg_KL$. \epf
\begin{center}
  {\bf  Acknowledgement}
  \end{center}
The author is very thankful to R.Y. Sharp for sending reprints of
his articles.
   \begin{center}
  {\bf REFERENCES}
  \end{center}
  \begin{enumerate} \item M.F.Atiyah, I.G.
Macdonald, Introduction to Commutative Algebra,Addison-Wesley
Publ.Co., 1969.
\item David Eisenbud, Commutative Algebra with a view Toward
Algebraic Geometry, Springer-Verlag, New York, Inc.,1995.
\item Arno van den Essen, Polynomial automorphims and the Jacobian
conjecture, Progress in Mathematics, Vol 190, Birkh\"{a}user, 2000.
\item Irving Kaplansky, Commutative Rings, The  University of
Chicago Press, Chicago, 1974.
\item Hideyuki Matsumura, Commutative Ring Theory, Cambridge
University Press, 1986.
\item Rodney Y. Sharp, Dimension of the
tensor product of two field extensions, Bulletin London Math. Soc.
9(1977), 42-48.
\item Rodney Y. Sharp and Peter Vamos, The dimension of the tensor
product of a finite number of field extensions., Jour. of Pure and
Applied Algebra, 10(1977), 249-252.
\end{enumerate}
\end{document}